\date{April 30, 2014}

\documentclass[a4paper,10pt,reqno]{amsart}
\usepackage{latexsym,amsmath,amsfonts,amscd,amssymb}
\usepackage{graphicx}
\usepackage{eucal}
\usepackage[all]{xy}
\usepackage{booktabs}

\newcommand\bC{\mathbb{C}}
\newcommand\bP{\mathbb{P}}
\newcommand\bQ{\mathbb{Q}}
\newcommand\bR{\mathbb{R}}
\newcommand\bZ{\mathbb{Z}}
\newcommand\cA{\mathcal{A}}
\newcommand\fg{\mathfrak{g}}

\newcommand\fk{\mathfrak{k}}
\newcommand\ba{\backslash}
\newcommand\ox{\otimes}
\newcommand\bk{\mathbf{k}}
\newcommand\id{\mathrm{Id}}
\newcommand\x{\times}
\newcommand\CP{\bC P}
\newcommand\imat{\mathbf{i}}
\renewcommand\Re{\mathrm{Re}}
\renewcommand\Im{\mathrm{Im}}

\newtheorem{theorem}{Theorem}
\newtheorem{lemma}[theorem]{Lemma}
\newtheorem{proposition}[theorem]{Proposition}

\newtheorem{corollary}[theorem]{Corollary}
\newtheorem{remark}[theorem]{Remark}
\theoremstyle{definition}

\newcommand{\thmref}[1]{Theorem~\ref{#1}}
\newcommand{\propref}[1]{Proposition~\ref{#1}}

\author[G. Bazzoni]{Giovanni Bazzoni}
\address{Fakult\"{a}t f\"{ur} Mathematik, Universit\"{a}t Bielefeld, Postfach 100301, D-33501 Bielefeld}

\email{gbazzoni@math.uni-bielefeld.de}

\author[V. Mu\~{n}oz]{Vicente Mu\~{n}oz}
\address{Facultad de Matem\'aticas, Universidad Complutense de Madrid, Plaza de Ciencias
3, 28040 Madrid, Spain}

\email{vicente.munoz@mat.ucm.es}


\title[Manifolds which are complex and symplectic but not K\"ahler]{Manifolds which are complex and symplectic but not K\"ahler}

\begin{document}

\begin{abstract}
The f{}irst example of a compact manifold admitting both complex and symplectic structures
but not admitting a K\"ahler structure is the renowned Kodaira-Thurston manifold. We review
its construction and show that this paradigm is very general and is not related to the
fundamental group. More specif{}ically, we prove that the simply-connected 
$8$-dimensional compact manifold of \cite{FM} admits both symplectic and complex structures
but does not carry K\"ahler metrics.
\end{abstract}

\maketitle

\section{Introduction} \label{sec:intro}

A complex manifold $M$ is a topological space modeled on open subsets of $\bC^n$ and with change
of charts being complex-dif{}ferentiable (that is, biholomorphisms). Here we say that $n$ is
the complex dimension of $M$. Complex manifolds are the objects that appear naturally in 
Algebraic Geometry: a projective variety is the zero locus of a collection of polynomials in the complex projective space $\CP^N$. 
When a projective variety is smooth and of complex dimension $n$, it is a complex manifold of dimension $n$.

A complex manifold $M$ of complex dimension $n$ is in particular a smooth dif{}ferentiable manifold of 
real dimension $2n$. Multiplication by $\imat$ on each complex tangent space $T_pM$, $p\in M$, gives
an endomorphism $J\colon TM\to TM$ such that $J^2=-\id$. An endomorphism $J\colon TM\to TM$ with
$J^2=-\id$ is called an \emph{almost complex} structure.  For a complex manifold $M$, $J$ satisf{}ies
that the Nijenhuis tensor 
  \begin{equation*} 
 N_J(X,Y)=[X,Y]+J[JX,Y]+J[X,JY]-[JX,JY]
 \end{equation*}
vanishes, $N_J(X,Y)=0$ for all vector f{}ields $X,Y$. In this case, we say that the almost complex structure
is \emph{integrable}. The celebrated Newlander-Nirenberg theorem \cite{NN} says that an almost
complex structure with $N_J=0$ is equivalent to a complex structure. Hence, for a smooth manifold $M$ 
to admit a complex structure, we need to check if there exist almost complex structures (this is a 
topological question), and then to f{}ind an integrable one (this is an analytic problem).

Projective varieties have further geometric properties. The complex proyective space $\CP^N$ has a 
natural hermitian metric, the Fubini-Study metric. This is the natural metric when one views $\CP^N$ as the homogeneous space
$U(N+1)/U(1)\x U(N)$. Therefore a projective variety $M\subset \CP^N$ inherits this hermitian
metric. Denote by $h$ the hermitian metric on $M$ and write $h=g+\imat \, \omega$, where
$g(X,Y)=\Re (h(X,Y))$ and $\omega(X,Y)=\Im (h(X,Y))=\Re(-\imat h(X,Y))=\Re (h(JX,Y))=g(JX,Y)$. 
Then $g$ is a Riemannian metric for which $J$ is an isometry ($g(JX,JY)=g(X,Y)$) and $\omega$ turns
out to be skew-symmetric, hence it is a $2$-form with $\omega(JX,JY)=\omega(X,Y)$ and $g(X,Y)=\omega(X,JY)$.
We say that $\omega$ is the fundamental form of $(M,h)$. This $2$-form is positive, in the sense
that $\omega^n>0$ (it gives the natural complex orientation). The Fubini-Study metric $h_{FS}$ has
fundamental form $\omega_{FS}\in \Omega^2(\CP^N)$. It is easy to see, using the $U(N+1)$-invariance,
that $d\omega_{FS}=0$. Therefore, for $\omega=\omega_{FS}|_M$ it also holds $d\omega=0$. 

We say that $(M,h)$ is a K\"ahler manifold when $M$ is a complex manifold and the fundamental form
$\omega$ satisf{}ies $d\omega=0$. A smooth projective variety is a K\"ahler manifold. Actually the converse 
holds when $[\omega]\in H^2(M,\bR)$ is an integral cohomology class, by Kodaira's theorem \cite{Wells}.

A dif{}ferent weakening of the K\"ahler condition (forgetting $J$ but keeping $\omega$) is that of a symplectic structure. A symplectic structure on 
a smooth $2n$-dimensional manifold $M$ is given by a $2$-form $\omega\in\Omega^2(M)$ which is 
closed ($d\omega=0$) and non-degenerate ($\omega^n$ is nowhere zero). Let $M$ be an even-dimensional 
manifold endowed with a complex structure $J$ and a symplectic structure $\omega$. Then $J$ is said to be 
\emph{compatible} with $\omega$ if, for vector f{}ields $X,Y$ on $M$, the bilinear form
\begin{equation}\label{compatibility}
  g(X,Y)=\omega(X,JY)   
\end{equation}
is a Riemannian metric. Therefore a K\"ahler manifold is a symplectic manifold endowed with a compatible 
complex structure, and $h=g+\imat\, \omega$ is the K\"ahler metric. The existence of a K\"ahler metric on a compact manifold constraints the
topology. In particular, if $(M,J,\omega)$ is a compact K\"ahler manifold of dimension $2n$, then (see \cite{Amoros,DGMS,GH,Wells})
\begin{enumerate}
 \item the fundamental group $\pi_1(M)$ belongs to a very restricted class of groups, called \emph{K\"ahler groups};
 \item $b_{2i-1}(M)$ is even for $i=1,\ldots,n$;
 \item the Lefschetz map $\mathcal{L}^{n-p}\colon H^p(M;\bR)\to H^{2n-p}(M;\bR)$, $a\mapsto [\omega]^{n-p}\wedge a$, is an isomorphism;
 \item $M$ is formal in the sense of Sullivan (see Section \ref{Section:Formality} for details).
\end{enumerate}

So it is natural to ask if the classes of smooth manifolds admitting complex, symplectic and K\"ahler structures coincide 
under some topological constraints. 

The lack of examples in symplectic geometry has been haunting this area of mathematics for many years now (pretty much since its \emph{d\'ebut} as a discipline in
its own). Indeed, the main source of examples of symplectic manifolds is Algebraic Geometry. This led to the belief that symplectic and K\"ahler conditions coincided in the compact case (see for
instance \cite{Gug}). There was a discrete breakthrough, in 1976, when Thurston \cite{Thu} gave the f{}irst example of a compact symplectic manifold with no K\"ahler structure.
Thurston's example had already been discovered, as a complex manifold, by Kodaira during his work on the classif{}ication of compact complex surfaces \cite{Kod}. We call it the
\emph{Kodaira-Thurston} manifold $KT$. Since $KT$ is a compact complex and symplectic manifold without K\"ahler structure, we obtain

\begin{theorem}\label{C_S_N_K}
 There exist compact manifolds which admit complex and symplectic structures but carry no K\"ahler metrics.
\end{theorem}

\begin{figure} \label{figu}
\begin{center}
\includegraphics[width=6cm]{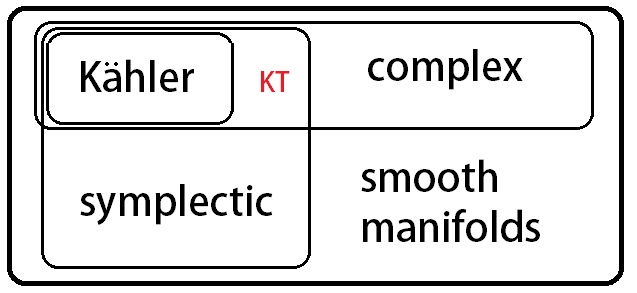}
\caption{Diagram of the dif{}ferent classes of manifolds, including $KT$}
\end{center}
\end{figure}

This means that the complex and symplectic structures that $KT$ admits cannot be compatible. The manifold $KT$ is in the place shown in Figure \ref{figu}.

The next natural question is whether some topological constraints may force the symplectic category to reduce to the K\"ahler one.
Regarding the fundamental group, it is natural to look for simply connected symplectic compact manifolds.
In \cite{McD}, McDuf{}f constructed a compact, simply connected, symplectic manifold with $b_3=3$, hence not K\"ahler.
For a detailed study on the relationship between formality and Lefschetz property on symplectic
manifolds, we refer to \cite{Cavalc}. In \cite{Bock}, Bock constructed non formal symplectic manifolds with arbitrary Betti numbers.

The construction of simply connected symplectic non formal (compact) manifolds turned out to be a more dif{}f{}icult problem. In fact,  it was conjectured in 1994 (see \cite{LO}) that a compact simply
connected symplectic (compact) manifold should be formal: this is the so-called \textit{Lupton-Oprea conjecture} on the formalising tendency of a symplectic structure. This conjecture was proven false
by Babenko and Ta\v{\i}manov in 2000 (see \cite{BT}). For every $n\geq 5$, they constructed an example of a simply connected, symplectic non formal compact manifold of real dimension $2n$.
On the other hand, by a result of Miller \cite{FM2,Mill}, simply connected compact manifolds of dimension $\leq 6$ are formal. 
Hence a remarkable gap in dimension $8$ was left. This gap was f{}illed by M. Fern\'andez and the second author in 2008 (see \cite{FM}).

Here we shall prove that the manifold constructed in \cite{FM} admits a complex structure, thereby giving a new example 
f{}itting in the scheme of \thmref{C_S_N_K}. The precise result is:

\begin{theorem}\label{main}
There exists an $8$-dimensional, compact, simply connected, symplectic and complex manifold 
which is non-formal and does not satisfy the Lefschetz property. In particular, it does not admit K\"ahler structures.
\end{theorem}


This paper is organized as follows. In Section \ref{Section:Formality} we recall the basics of rational homotopy theory and formality. 
In Section \ref{Section:Kodaira-Thurston} we give a description of $KT$, construct explicit complex and symplectic structures on it 
and show that it carries no K\"ahler metric. In Section \ref{Section:Review} we review the construction of 
the symplectic manifold $(\widetilde{M},\widetilde{\omega})$ of \cite{FM}. This  is constructed by resolving symplectically the singularities of a symplectic orbifold $(\widehat{M},\widehat{\omega})$,
a quotient of a compact symplectic nilmanifold $(M,\omega)$ by a certain $\bZ_3$-action. In Section \ref{Section:Complex_structure} we describe a complex structure $\widehat{J}$ on the orbifold
$\widehat{M}$ and construct a complex resolution of singularities $(\overline{M},\overline{J})$. Finally, in Section \ref{Section:Diffeomorphism} we show that $\widetilde{M}$ and $\overline{M}$ are
dif{}feomorphic.


{\bf Acknowledgements} The main result of this paper has been taken 
from Chapter 1 of the PhD thesis \cite{thesis} of the f{}irst author.
We are grateful to Jes\'us Ruiz who suggested us this question. 
Partially supported by (Spain) MICINN grant MTM2010-17389.

\section{Formality}\label{Section:Formality}

Formality is a property of the rational homotopy type of a space $X$. We present here a rough introduction, 
referring to \cite{FHT,FOT,GM} for more details. By space, we mean a connected CW
complex of f{}inite type (we allow a f{}inite number of cells in each dimension) which is nilpotent 
(its fundamental group is nilpotent and acts nilpotently on higher homotopy groups). A space $X$ is
rational if $\pi_i(X)$ is a rational vector space for every $i\geq 1$ (recall that a nilpotent group has 
a well def{}ined rationalization). The rationalization of a space $X$ is a rational space $X_\bQ$
together with a map $f\colon X\to X_\bQ$ such that $f_i\colon\pi_i(X)\ox\bQ\to\pi_i(X_\bQ)$ is an 
isomorphism for every $i\geq 1$. We identify two spaces if they have a common rationalization.
By rational homotopy type of a space $X$ we mean the homotopy type of its rationalization. 
Quillen and Sullivan proposed two dif{}ferent approaches to capture the rational homotopy type of a
space in an algebraic model, see \cite{Q,Su}. Here we review briefly Sullivan's ideas.

A \emph{commutative dif{}ferential graded algebra} $(\cA,d)$ over a f{}ield $\bk$ of zero characteristic ($\bk$-cdga for short) 
is a graded algebra $\cA=\oplus_{i\geq 0}A^i$ which is graded commutative,
together with a $\bk$-linear map $d\colon A^i\to A^{i+1}$, the dif{}ferential, which satisf{}ies $d^2=0$ and which is a 
graded derivation, i.e.,\ for homogeneous elements $a\in A^p$ and $b\in A^q$,
\[
 d(a\cdot b)=(da)\cdot b+(-1)^{pq}a\cdot(db).
\]
The cohomology of a $(\cA,d)$, denoted $H^*(\cA)$, is a $\bk$-cdga with trivial dif{}ferential. A $\bk$-cdga is \emph{connected} if $H^0(\cA)\cong\bk$.

The de Rham algebra $\Omega(M)$ of a smooth manifold $M$, together with the exterior dif{}ferential, is an $\bR$-cdga. 
The piecewise linear forms $A_{PL}(X)$ on a PL-manifold $X$, endowed
with a suitable dif{}ferential combining the exterior dif{}ferential and the boundary of simplices, form a $\bQ$-cdga (see \cite{GM}). 
There is a de Rham-type theorem for both cdga's, hence we have isomorphisms
\[
 H^*(\Omega(M))\cong H^*(M;\bR) \quad \mathrm{and} \quad H^*(A_{PL}(X))\cong H^*(X;\bQ).
\]

Let $X$ be a space. The idea of Sullivan is to replace $A_{PL}(X)$ by another $\bQ$-cdga, which has 
the same cohomological information as $A_{PL}(X)$ but is algebraically more tractable: the
\emph{minimal model}. A $\bk$-cdga $(\cA,d)$ is \emph{minimal} if
\begin{itemize}
 \item $\cA$ is the free graded algebra over a graded vector space $V=\oplus_iV^i$; this means that 
$\cA$ is the tensor product of the exterior algebra on the odd degree generators and the
 symmetric algebra on the even degree generators, $\cA=\mathrm{Ext}(V^{\mathrm{odd}})\ox \mathrm{Sym}(V^{\mathrm{even}})$. 
The standard notation is $\cA=\Lambda V$.
 \item there exists a collection $\{x_i\}_{i\in\mathcal{I}}$ of generators of $V$, indexed by a well-ordered set 
$\mathcal{I}$, such that $|x_i|\leq|x_j|$ if $i<j$ and the dif{}ferential of a generator
 $x_j$ is an element of $\Lambda(V^{<j})$. Here $|\cdot|$ denotes the degree and $V^{<j}$ consists of 
the generators $x_i$ with $i<j$. Notice, in particular, that $d$ does not have linear part.
\end{itemize}

We denote a minimal $\bk$-cdga by $(\Lambda V,d)$. A \emph{minimal model} for a $\bk$-cdga $(\cA,d)$ is a 
minimal $\bk$-cdga $(\Lambda V,d)$ together with a $\bk$-cdga morphism $\phi\colon
(\Lambda V,d)\to  (\cA,d)$ which induces an isomorphism in cohomology (such a morphism is called \emph{quasi-isomorphism}).

We have the following fundamental result:

\begin{theorem}[\cite{FHT}, Theorem 14.12]
 Any connected $\bk$-cdga has a minimal model, which is unique up to isomorphism.
\end{theorem}

By def{}inition, the rational minimal model of a space $X$, $(\Lambda V_X,d)$, is the minimal 
model of the $\bQ$-cdga $(A_{PL}(X),d)$. One can show that, when $M$ is a smooth
manifold, the real minimal model of $M$ can be computed from the de Rham algebra 
$(\Omega(M),d)$. A central result in rational homotopy theory is the following:

\begin{theorem}[\cite{Su}]\label{fundamental}
 Two spaces have the same rational homotopy type if and only if their rational minimal models are isomorphic.
\end{theorem}

In particular, PL forms (resp.\ smooth forms) contain all the rational-homotopic (resp.\ real-homotopic) information of a space (smooth manifold). It is often dif{}f{}icult to know
the whole de Rham algebra of a manifold; it would be very convenient if the (say, real) minimal model could be constructed directly from the de Rham cohomology. A space for which this happens is
called \emph{formal}. More precisely, a space $X$ is formal if there exists a quasi-isomorphism $(\Lambda V_X,d)\to (H^*(X;\bQ),0)$. In particular, the rational homotopy type of a formal space $X$ is
a formal consequence of its rational cohomology. Many spaces are known to be formal: compact Lie groups, $H$-spaces, symmetric spaces, $\ldots$ 
For us, the relevant result is the following:

\begin{theorem}[\cite{DGMS}]
  A smooth compact manifold $M$ admitting a K\"ahler structure is formal.
\end{theorem}

A very useful criterion for establishing formality is the following:

\begin{theorem}[\cite{DGMS}, Theorem 4.1]\label{C+N}
 Let $X$ be a space and let $(\Lambda V_X,d)$ be its minimal model. Then $X$ is formal if and only if we can write 
$V_X=C\oplus N$ with $d=0$ on $C$ and $d$ injective on $N$, in such way that every closed
element in the ideal generated by $N$ in exact.
\end{theorem}

Let $(\cA,d)$ be a $\bk$-cdga and let $H^*(\cA)$ be its cohomology. Let $a\in H^{|a|}(\cA)$, $b\in H^{|b|}(\cA)$ and
$c\in H^{|c|}(\cA)$ such that $a\cdot b=b\cdot c=0$. Then $a\cdot b\cdot c$ is zero for two reasons. 
Consequently, a dif{}ference element $\langle a,b,c\rangle\in H^{|a|+|b|+|c|-1}(\cA)/\mathcal{J}$ can be formed,
where $\mathcal{J}$ is the ideal generated by $a$ and $c$ in $H^*(\cA)$. Take cocycles 
$\alpha,\beta,\gamma\in \cA$ representing $a,b,c$ respectively. Then $\alpha\cdot\beta=d\xi$ and
$\beta\cdot\gamma=d\eta$, hence $\xi\cdot\gamma+(-1)^{|a|+1}\alpha\cdot\eta$ is a 
closed $(|a|+|b|+|c|-1)$-form whose cohomology class is well def{}ined modulo $\mathcal{J}$. We set 
$\langle a,b,c\rangle=[\xi\cdot\gamma+(-1)^{|a|+1}\alpha\cdot\eta]$. Then $\langle a,b,c\rangle$ is called
the \emph{triple Massey product} of the cohomology classes $a,b,c$. 

The def{}inition of higher Massey products is as follows (see~\cite{K,Mas}). Given $a_i\in H^{|a_i|}(\cA)$, $1\leq i\leq t$, $t\geq 3$, the Massey product $\langle
a_1,a_2,\dots,a_t\rangle$ is def{}ined if there are $\alpha_{i,j}\in \cA$, 
with $1\leq i\leq j\leq t$, except for the case $(i,j)=(1,t)$, such that
 \begin{equation}\label{eqn:gm}
 a_i=[\alpha_{i,i}], \qquad d\alpha_{i,j}= \sum\limits_{k=i}^{j-1} (-1)^{|\alpha_{i,k}|}  {\alpha}_{i,k}\wedge \alpha_{k+1,j}.
 \end{equation}
Then the Massey product is
 \begin{equation}\label{eqn:gm2}
 \langle a_1,a_2,\dots,a_t \rangle =\left\{
 \left[\sum\limits_{k=1}^{t-1} (-1)^{|\alpha_{1,k}|}{\alpha}_{1,k} \wedge
 \alpha_{k+1,t}\right]\right\}
 \subset H^{|a_1|+ \cdots +|a_t|-(t-2)}(\cA)\, ,
 \end{equation}
where the $\alpha_{i,j}$ are as in \eqref{eqn:gm}. We say that the Massey product is trivial if $0\in \langle a_1,a_2,\dots,a_t\rangle$. Note that for $\langle a_1,a_2,\dots,a_t\rangle$ to
be def{}ined it is necessary that both $\langle a_1,\dots,a_{t-1}\rangle$ and $\langle a_2,\dots,a_t\rangle$ are def{}ined and trivial.

\begin{proposition} \label{prop:Massey}
 If $X$ is formal then all (higher) Massey products of $(\Lambda V_X,d)$ are zero.
\end{proposition}

\begin{proof}
 The proof can be found in \cite{BT2}. We shall give a simple proof for the case of
triple and quadruple Massey products, which suf{}f{}ices for this paper.

 As $X$ is formal, \thmref{C+N} guarantees that we can write $V_X=C\oplus N$ with $d=0$ on $C$ and $d$ injective on $N$, 
 in such way that every closed element in the ideal $I(N)$ generated by $N$ in exact.
 Note that there is a decomposition $\Lambda V=\Lambda C \oplus I(N)$. 
 Let $a_i\in H^{|a_i|}(\cA)$, $1\leq i\leq t$. By def{}inition of Massey product, there are $\alpha_{i,i}\in \Lambda V$ 
 with $a_i=[\alpha_{i,i}]$, and for each $i<j$, $(i,j)\neq (1,t)$, there are $\alpha_{i,j}$ with 
 $d\alpha_{i,j}= \sum\limits_{k=i}^{j-1} (-1)^{|\alpha_{i,k}|}  {\alpha}_{i,k}\wedge \alpha_{k+1,j}$.

Write $\alpha_{i,j}=\beta_{i,j}+\eta_{i,j}$ with $\beta_{i,j}\in \Lambda C$, $\eta_{i,j}\in I(N)$. 
As $d\alpha_{i,j}=d\eta_{i,j}$, we can use in the case of triple Massey products (that is, $t=3$),
the elements $\eta_{12}$ and $\eta_{23}$. Then the triple Massey product $\langle a_1,a_2,a_3\rangle$ contains
$ (-1)^{|\eta_{12}|} \eta_{12}\alpha_{33}+ (-1)^{|\alpha_{11}|} \alpha_{11}\eta_{23}$ which 
is in $I(N)$, hence exact.

In the case of quadruple Massey products (that is, $t=4$), we use $\eta_{12}, \eta_{23},\eta_{34}$ instead of
$\alpha_{12},\alpha_{23},\alpha_{34}$. The equation 
 \begin{align*}
 d\alpha_{13} &= 
(-1)^{|\alpha_{12}|} \alpha_{12}\alpha_{33}+ (-1)^{|\alpha_{11}|} \alpha_{11}\alpha_{23} \\
 &=
(-1)^{|\alpha_{12}|} \eta_{12}\alpha_{33}+ (-1)^{|\alpha_{11}|} \alpha_{11}\eta_{23}+
(-1)^{|\alpha_{12}|} \beta_{12}\alpha_{33}+ (-1)^{|\alpha_{11}|} \alpha_{11}\beta_{23}
 \end{align*}
implies
that $(-1)^{|\alpha_{12}|} \eta_{12}\alpha_{33}+ (-1)^{|\alpha_{11}|} \alpha_{11}\eta_{23}$ is
closed, hence exact (as it lives in $I(N)$). Write it as $d\psi_{13}$ with $\psi_{13} \in I(N)$. Analogously
def{}ine $\psi_{24}\in I(N)$. Thus  the quadruple Massey product  $\langle a_1,a_2,a_3,a_4\rangle$ contains
$ (-1)^{|\psi_{13}|} \psi_{13}\alpha_{44}+ (-1)^{|\eta_{12}|} \eta_{12}\eta_{34}+ (-1)^{|\alpha_{11}|} \alpha_{11}\psi_{24}$ 
which is in $I(N)$, hence exact.
 \end{proof}

\section{The Kodaira-Thurston manifold}\label{Section:Kodaira-Thurston}

The Kodaira-Thurston manifold can be described in various ways. For Kodaira, $KT$ was a compact quotient of $\bC^2$ by a certain group acting co-compactly. 
Thurston interpreted it as a symplectic $T^2$-bundle over
$T^2$. In this section we describe it as a nilmanifold, write down explicit symplectic and complex structures on $KT$ and show that $KT$ carries no K\"ahler metric. 

A nilmanifold is a compact quotient of a simply connected, nilpotent Lie group $G$ by a lattice $\Gamma$. 
Since $\Gamma$ is a subgroup of a nilpotent group, it is also nilpotent. The exponential map
$\exp\colon\fg\to G$ is a dif{}feomorphism, hence $G\cong\bR^n$ for some $n$. Therefore, if $N=\Gamma\ba G$ is a compact nilmanifold, $G\to N$ is the universal cover, $\pi_1(N)\cong\Gamma$ and
$\pi_i(N)=0$
for $i\geq 2$. Hence a nilmanifold is a nilpotent space.

Nilmanifolds are interesting because they are a rich source of answers to many questions in dif{}ferent areas of Mathematics. As we already mentioned, $KT$ was the f{}irst example of a compact
symplectic
non K\"ahler manifold. From the point of view of complex geometry, there exist complex nilmanifolds for which the Fr\"olicher spectral sequence is arbitrarily
non-degenerate, see \cite{Rol}.

K\"ahler nilmanifolds are very special:

\begin{theorem}[Benson-Gordon, Hasegawa \cite{BG,Has}]\label{BG}
 Let $N$ be a compact symplectic nilmanifold endowed with a K\"ahler structure. Then $N$ is dif{}feomorphic to a torus.
\end{theorem}

Benson and Gordon proved that a symplectic nilmanifold $N$ of dimension $2n$ for which the Lefschetz map 
$\mathcal{L}^{n-1}\colon H^1(N;\bR)\to H^{2n-1}(N;\bR)$ is an isomorphism is dif{}feomorphic to a torus. Hasegawa showed
that a formal nilmanifold $N$ is dif{}feomorphic to a torus. Notice, however, that there exists 
many examples of non-toral symplectic and complex nilmanifolds (see \cite{BM,Gu,Sal}).

Let $H$ denote the Heisenberg group, i.e.
\[
H=\left\{\begin{pmatrix}
1 & b & c\\
0 & 1 & a\\
0 & 0 & 1
\end{pmatrix} \ | \ a,b,c \in\bR\right\}
\]
and let $H_{\bZ}$ denote the subgroup of matrices with entries in $\bZ$. Then $H$ is a nilpotent Lie group, 
dif{}feomorphic to $\bR^3$, $H_{\bZ}\subset H$ is a lattice and $N=H_\bZ\backslash H$ is a
compact nilmanifold. Let $G=H\times\bR$ and $G_\bZ=H_\bZ\times\bZ$. The Kodaira-Thurston manifold is $KT=G_\bZ\backslash G$.

Let $\fk$ be a Lie algebra over a f{}ield $\bk$ of characteristic zero. 
The exterior algebra $\Lambda \fk^*$ is endowed with a dif{}ferential $d\colon\Lambda^p\fk^*\to\Lambda^{p+1}\fk^*$, def{}ined as
follows: $d\colon\fk^*\to\Lambda^2\fk^*$ is the dualization of the bracket, i.e. $(d\alpha)(X,Y)=-\alpha([X,Y])$ if $\alpha\in\fk^*$ 
and $X,Y\in\fk$. $d$ is then extended to $\Lambda \fk^*$ by
imposing the graded Leibnitz rule: for $\alpha\in\Lambda^p\fk^*$ and $\beta\in\Lambda^q\fk^*$, 
$d(\alpha\wedge\beta)=(d\alpha)\wedge\beta+(-1)^{pq}\alpha\wedge(d\beta)$. The vanishing of $d^2$ is
equivalent to the Jacobi identity in $\fk$. In the language of the Section \ref{Section:Formality}, $(\Lambda\fk^*,d)$ 
is a $\bk$-cdga, known as Chevalley-Eilenberg complex of $\fk$.

 
Let $\fg$ be the Lie algebra of $G$ and let $\fg^*$ be its dual. We identify tensors on $\fg$ and 
$\fg^*$ with left-invariant objects on $G$, which therefore descend to $KT$. It is easy to check that
$\fg$ has a basis $\langle X_1,X_2,X_3,X_4\rangle$ in which the only non-zero bracket is 
$[X_1,X_2]=-X_3$. Let $\langle x_1,x_2,x_3,x_4\rangle$ be the dual basis of $\fg^*$. The
only non-zero dif{}ferential on $\fg^*$ is computed to be $dx_3=x_1\wedge x_2$. 

The element $\omega=x_1\wedge x_4+x_2\wedge x_3\in\Lambda^2\fg^*$ is closed and non-degenerate. By abuse of notation, we
denote by $\omega$ the corresponding left-invariant symplectic structure on $KT$ as well. Thus $(KT,\omega)$ is a compact symplectic 4-manifold.

Recall that if $\fk$ is an even-dimensional Lie algebra, $J\colon\fk\to\fk$ is a complex structure if $J^2=-\id$ and it satisf{}ies the integrability condition
\begin{equation}\label{Nijenhuis}
N_J(X,Y)=[X,Y]+J[JX,Y]+J[X,JY]-[JX,JY]=0, \ \mathrm{for} \ X,Y\in\fk.
\end{equation}

In our situation, def{}ine $J\colon\fg\to\fg$ by
\[
J(X_1)=X_2, \ J(X_2)=-X_1, \ J(X_3)=X_4 \ \mathrm{and} \ J(X_4)=-X_3.
\]
A straightforward computation shows that \eqref{Nijenhuis} holds, hence $J$ is a complex 
structure on $\fg$. Again by abuse of notation, we denote by $J$ the corresponding
left-invariant complex structure on $KT$. Thus $(KT,J)$ is a compact complex surface.

Let $N=\Gamma\ba G$ be a compact nilmanifold. Considering $\Lambda\fg^*$ as left-invariant 
forms on $N$, we have a natural inclusion $\iota\colon(\Lambda\fg^*,d)\to(\Omega(N),d)$. By a result of
Nomizu (see \cite{Nom}), $\iota$ is a quasi-isomorphism, hence the de Rham cohomology of $N$ is isomorphic to the cohomology of the Chevalley-Eilenberg complex of $\fg$. In our case, three of the four
generators of $\fg^*$ are closed, hence we get $b_1(KT)=3$.

Since $KT$ has an odd Betti number which is odd, we see that it does not carry any K\"ahler metric. 
We also see explicitly that $KT$ does not satisfy the Lefschetz property. Indeed, take $[x_2]\in
H^1(KT;\bR)$. Then $\mathcal{L}\colon H^1(KT;\bR)\to H^3(KT;\bR)$ sends $[x_2]$ to $-[x_1\wedge x_2\wedge x_4]=-[d(x_3\wedge x_4)]=0$.

The Lie algebra $\fg$ is endowed with a complex structure $J$ and a symplectic structure $\omega$. 
Def{}ine a tensor $g\colon\fg\ox\fg\to\bR$ by
\[
g(X,Y)=\omega(X,JY), \ X,Y\in\fg.
\]
It is easy to see that the matrix of $g$ in the basis $\langle X_1,X_2,X_3,X_4\rangle$ is
\[
\begin{pmatrix}
0 & 0 & 1 & 0\\
0 & 0 & 0 & -1\\
-1 & 0 & 0 & 0\\
0 & 1 & 0 & 0\\
\end{pmatrix}.
\]
$g$ is not a scalar product on $\fg$, hence the corresponding left-invariant tensor on $KT$ is not a Riemannian metric.

Let $M$ be a manifold endowed with a complex structure $J$ and a symplectic structure $\omega$. One could in principle relax condition \eqref{compatibility} above and ask $J$ to be only \emph{tamed}
by $\omega$, which means $\omega(X,JX)>0$ for $X\in\mathfrak{X}(M)$. A symplectic manifold $(M,\omega)$ endowed with a tamed complex structure $J$ is called \emph{Hermitian-symplectic}. There are no
known examples of compact Hermitian-symplectic non K\"ahler manifolds.

We see that $(KT,J,\omega)$ is not Hermitian-symplectic. Indeed, $\omega(X_1,JX_1)=0$. It is proved in \cite{EFV} that a compact nilmanifold endowed with a Hermitian-symplectic structure
is actually K\"ahler. Hence we see that $KT$ does not carry \emph{any} Hermitian-symplectic structure (not just left-invariant).

To see explicitly that $KT$ is non formal, we need to compute the minimal model of a nilmanifold.

\begin{theorem}[\cite{Has}]\label{Hasegawa}
 Let $N=\Gamma\ba G$ be a compact nilmanifold. Then $(\Lambda \fg^*,d)$ is the rational minimal model of $N$.
\end{theorem}

Since a nilmanifold is a nilpotent space, Theorem \ref{fundamental} holds and the rational homotopy of a compact nilmanifold is codif{}ied in the corresponding minimal model. 
Here $(\Lambda \fg^*,d)$ is  a minimal algebra generated in degree $1$. By a result of Mal'cev (see \cite[Theorem 2.12]{Rag}), 
a simply connected nilpotent Lie group $G$ admits a lattice if
and only if $\fg$ admits a basis such that the structure constants are rational numbers. 
Hence, if $N=\Gamma\ba G$ is a compact nilmanifold, $(\Lambda \fg^*,d)$ is automatically a $\bQ$-cdga.

Applying Theorem \ref{Hasegawa}, the minimal model of $KT$ is
\[
 (\Lambda^*\langle x_1,x_2,x_3,x_4\rangle, dx_3=x_1\wedge x_2).
\]
In the notation of Theorem \ref{C+N}, we have $C=\langle x_1,x_2,x_4\rangle$ and $N=\langle x_3\rangle$. The element $x_1\wedge x_3$ belongs to the ideal generated by $N$, is closed, but
not exact. A nonzero Massey product is constructed as follows. Take $a=b=[x_1]$ and $c=[x_2]$ in $H^1(KT;\bQ)$. The recipe given after Theorem \ref{C+N} tells us that the triple Massey product
$\langle [x_1],[x_1],[x_2]\rangle=[x_1\wedge x_3]$ is a well def{}ined element of $H^2(KT;\bQ)$, which is non-zero modulo the ideal generated in $H^*(KT;\bQ)$ by $[x_1]$ and $[x_2]$.

\section{A simply connected symplectic non-formal 8-manifold}\label{Section:Review}

In this section we recall the construction of a simply connected, 8-dimensional symplectic non formal manifold performed in \cite{FM}. Although quite involved, this construction also starts with a
nilmanifold, 
showing the importance of such manifolds in the whole theory.

Let $H_{\bC}$ be the complex Heisenberg group, def{}ined as
\[
H_{\bC}=\left\{
A=\begin{pmatrix}
1 & u_2 & u_3\\
0 & 1 & u_1\\
0 & 0 & 1\\
\end{pmatrix} \quad | \quad u_1,u_2,u_3\in\bC\right\}.
\]
The map $H_{\bC}\to\bC^3$, $A\mapsto(u_1,u_2,u_3)$, gives a global system of holomorphic 
coordinates on $H_{\bC}$. Set $G=H_{\bC}\times \bC$, with global coordinates $(u_1,u_2,u_3,u_4)$. Let
$\Gamma\subset\bC$ be the lattice generated by $1$ and $\zeta=e^{2\pi \imat /3}$. 
Also, let $G_\Gamma\subset G$ be the discrete subgroup of matrices with entries in $\Gamma$. We let $G_\Gamma$ act on $G$
on the left and set $M=G_\Gamma\backslash G$. Then $M$ is a compact complex parallelizable 
nilmanifold. Notice that $M$ can be seen as a principal torus bundle
\[
T^2=\Gamma\ba\bC \hookrightarrow M\rightarrow T^6=(\Gamma\ba\bC)^3
\]
using the projection $(u_1,u_2,u_3,u_4)\mapsto (u_1,u_2,u_4)$. $M$ is a complex version of the Kodaira-Thurston manifold.

We interpret $\bZ_3$ as the group of cubic roots of unity and consider the right $\bZ_3$-action $\rho\colon\bZ_3\times G\to G$ given, in terms of a generator $\zeta=e^{2\pi \imat/3}$, by
\begin{equation}\label{Z3-action}
(\zeta,(u_1,u_2,u_3,u_4))\mapsto (\zeta u_1,\zeta u_2,\zeta^2 u_3,\zeta u_4). 
\end{equation}
This action preserves the group operation on $G$ and the lattice, hence descends to an action of $\bZ_3$ on $M$. Set $\widehat{M}=M/\bZ_3$. 
Then $\widehat{M}$ is not smooth,  it has $81$ isolated quotient singularities.

A basis for left-invariant $1$-forms on $G$ is given by
\[
\mu=d u_1, \quad \nu=d u_2, \quad \theta=d u_3-u_2 du_1 \quad \text{and} \quad \eta=du_4
\]
(over the complex numbers), with
\[
d\mu=d\nu=d\eta=0, \quad d\theta=\mu\wedge\nu.
\]
The action of $\bZ_3$ on left-invariant 1-forms is given by
\[
\rho^*\mu=\zeta\mu, \quad \rho^*\nu=\zeta\nu, \quad \rho^*\theta=\zeta^2\theta \quad \text{and} \quad \rho^*\eta=\zeta\eta.
\]
The $2$-form
\begin{equation}\label{simplform}
 \omega=\imat\mu\wedge\bar{\mu}+\nu\wedge\theta+\bar{\nu}\wedge\bar{\theta}+\imat\eta\wedge\bar{\eta}
\end{equation}
on $M$ satisf{}ies $\bar{\omega}=\omega$, so it is real. It is closed and satisf{}ies $\omega^4\neq 0$. Thus $\omega$ is a symplectic form. Notice also that
\[
\rho^*\omega=\zeta^3(\imat\mu\wedge\bar{\mu}+\nu\wedge\theta+\imat\eta\wedge\bar{\eta})+\zeta^{-3}\bar{\nu}\wedge\bar{\theta}=\omega,
\]
hence $\omega$ is $\bZ_3-$invariant and descends to a symplectic form $\widehat{\omega}$ on the quotient $\widehat{M}$. 
Therefore $(\widehat{M},\widehat{\omega})$ is a symplectic orbifold. In \cite{FM} a desingularization procedure for the symplectic orbifold is given, 
producing a symplectic manifold. 

\begin{proposition}[\cite{FM}, Propositions 2.1 and 2.3]\label{prop:21}
There exists a smooth compact simply connected symplectic manifold $(\widetilde{M},\widetilde{\omega})$ 
which is isomorphic to $(\widehat{M},\widehat{\omega})$ outside a small neighborhood of the
singular points.
\end{proposition}

In \cite{FM}, it is shown that $(\widetilde{M},\widetilde{\omega})$ is non formal, and also that it does not satisfy the Lefschetz
property (see Remark 3.3 in \cite{FM}). The non-formality is seen in \cite{FM} via a modif{}ication of the Massey product, named
$a$-Massey product, which are studied extensively in  \cite{CFM}. Here we shall see the non-formality of $\widehat{M}$ by
showing that there exists a non-zero quadruple Massey product. Transfering the Massey product from $\widehat{M}$ to the desingularization
$\widetilde{M}$ follows the arguments of \cite[Theorem 3.2]{FM} and it is quite standard.

The complex minimal model of $M$ is $\Lambda V_M=\Lambda (\mu,\nu,\theta,\eta, \bar\mu,\bar\nu,\bar\theta,\bar\eta)$ with $d\theta=\mu\wedge\nu$ and
$d\bar\theta=\bar\mu\wedge\bar\nu$. Our orbifold is $\widehat{M}=M/\bZ_3$, where $\bZ_3$ acts in the minimal model as
$(\mu,\nu,\theta,\eta)\mapsto (\zeta\mu,\zeta\nu,\zeta^2\theta,\zeta\eta)$. A model (that is, a $\bC$-cdga quasi-isomorphic to its
minimal model) for $\widehat{M}$ is given by $\cA=(\Lambda V_M)^{\bZ_3}$. Easily 
 \begin{align*}
 A^1 &= 0, \\
 A^2 &=\big( \langle \mu,\nu,\eta\rangle \wedge\langle \bar\mu,\bar\nu,\bar\eta\rangle \big) \oplus 
  \langle \mu\wedge \theta,\nu\wedge \theta,\eta\wedge \theta , \bar\mu\wedge \bar\theta,\bar\nu\wedge \bar\theta,\bar\eta\wedge \bar\theta , \theta\wedge \bar\theta\rangle ,\\
 A^3 &= \Lambda^3(\mu,\nu,\eta,\bar\theta) \oplus \Lambda^3(\bar\mu,\bar\nu,\bar\eta,\theta) .
\end{align*}
With this, one can check that $H^3(\cA)=0$. 

Take now $a_1=[\nu\wedge \bar\eta]$, $a_2=[\mu\wedge \bar\mu]$, $a_3=[\mu\wedge \bar\mu]$ and 
$a_4=[\eta\wedge\bar\nu]$. We shall compute $\langle a_1,a_2,a_3,a_4\rangle$ 
and check that it does not contain the zero element. 
A Massey product $b\in \langle a_1,a_2,a_3,a_4\rangle$ is computed according to formula \eqref{eqn:gm2}.
As $A^1=0$, it must be $\alpha_{11}=\nu\wedge \bar\eta$, $\alpha_{22}=\alpha_{33}
=\mu\wedge \bar\mu$ and $\alpha_{44}=\eta\wedge\bar\nu$. Then
 \begin{align*}
 \alpha_{12}&= -\theta\wedge \bar\mu\wedge \bar\eta + z_1 , \\
 \alpha_{13}&= \nu\wedge\bar\eta\wedge f_2-f_1\wedge\mu\wedge\bar\mu + w_1 , \\
 \alpha_{23} &= z_2 , \\
 \alpha_{24}&= \mu\wedge\bar\mu\wedge f_3-f_2\wedge \eta\wedge\bar\nu +w_2 ,\\
 \alpha_{34}&= -\bar\theta\wedge \mu\wedge\eta +z_3,
\end{align*}
where $z_1,z_2,z_3\in A^3$ are closed, hence exact, thus $z_i=df_i$, with $f_i\in A^2$,
and $w_1,w_2\in A^4$ are closed. A computation gives
\begin{align*}
  b  =& [\alpha_{11}\wedge \alpha_{24} -\alpha_{12}\wedge \alpha_{34} + \alpha_{13}\wedge\alpha_{44}]
 \\
  = & [\theta\wedge\bar\theta\wedge\mu\wedge\bar\mu\wedge\eta\wedge\bar\eta
+ w_1\wedge\eta\wedge\bar\nu + w_2\wedge\nu\wedge\bar\eta].
\end{align*}
To check that this is non-zero, we multiply by $[\nu\wedge\bar\nu]$. Then the terms with $w_1$
and $w_2$ cancel, so $b\wedge [\nu\wedge\bar\nu]\neq 0$, hence $b\neq 0$.

\section{The complex structure}\label{Section:Complex_structure}

In this section we describe the complex structure $J$ on $G$ in two equivalent ways, and we show that it descends to $M=G_\Gamma\ba G$ and also to the orbifold
$\widehat{M}=(G_\Gamma\ba G)/\bZ_3$. Then we construct a complex resolution of singularities, which will give a smooth complex $4$-fold $(\overline{M},\overline{J})$.

Let us consider the group $G=H_{\bC}\times\bC$ above. $G$ can be realized as a complex Lie subgroup of 
$\text{GL}(5,\bC)$ by sending the pair $(A,u_4)\in H_{\bC}\times\bC$ to the matrix
\[
\begin{pmatrix}
1 & u_2 & u_3 & 0 & 0\\
0 & 1 & u_1 & 0 & 0\\
0 & 0 & 1 & 0 & 0\\
0 & 0 & 0 & 1 & u_4\\
0 & 0 & 0 & 0 & 1\\
\end{pmatrix}.
\]
$\text{GL}(5,\bC)$ is an open subset of $\bC^{25}$, hence each tangent space $T_X\text{GL}(5,\bC)\cong\bC^{25}$, 
$X\in\text{GL}(5,\bC)$, inherits the standard complex structure of the ambient space,
which is the multiplication by $\imat=\sqrt{-1}$. As a complex submanifold of $\text{GL}(5,\bC)$, 
$G$ inherits the same complex structure on each tangent space. This means that the complex structure on $G$
is multiplication by $\imat$ on each tangent space $T_gG$, $g\in G$. The left translations $L_g\colon G\to G$, $h\mapsto gh$, 
are holomorphic maps, since they are written as polynomials in local
coordinates. This shows that $G$ is a complex parallelizable Lie group: the dif{}ferential of $L_g$ is complex linear and 
a parallelization is given by moving $T_eG$ around. Let $J$ denote the complex
structure on $G$ induced by the inclusion $G\hookrightarrow\text{GL}(5,\bC)$. The above considerations show that $J$ is left invariant.

Let us consider the tangent space $T_eG$, where $e\in G$ is the identity. There is an identif{}ication between the 
Lie algebra $\fg$ of $G$ and the vector space of left invariant holomorphic
vector f{}ields on $G$, endowed with the natural Lie bracket. The complex structure on $\fg$ is multiplication 
by $\imat$, and $\fg$ is a complex Lie algebra of dimension $4$, described as follows:
\[
\fg=\{\langle Z_1,Z_2,Z_3,Z_4\rangle \ | \ [Z_1,Z_2]= -Z_3\}.
\]
By identifying $\fg$ with $T_eG$, one has $T_gG=d_eL_g(\fg), \ \forall g\in G$. This shows again that the 
complex structure $J_g$ on $T_gG$ is multiplication by $\imat$, for every $g\in G$.

We go through the details of the construction of  left invariant complex structure on $G$. Let $J_e$ 
denote the complex structure (i.e.\ multiplication by $\imat$) on $\fg$ and let $g\in G$ be a point.
Def{}ine the complex structure 
$J_g\colon T_gG\to T_gG$ as
\[
J_g(X(g))=d_eL_g(\imat x),
\]
where $X$ is a left invariant vector f{}ield on $G$ and $x\in\fg$ is such that $d_eL_g(x)=X(g)$. 
This def{}ines $J$ as a smooth section of the bundle $\text{End}(TG)$. Let us show that $J^2=-\id$.
Indeed,
\[
J_g^2(X(g))=J_g(J_g(X(g)))=d_eL_g(\imat(\imat x))=-d_eL_g(x)=-X(g).
\]
\begin{lemma}
The (almost) complex structure def{}ined above is left invariant.
\end{lemma}
\begin{proof} We must prove that, for every $g\in G$, $(L_g)^*J=J$. So take $X(h)\in T_h G$. Then
\[
J_h(X(h))=d_eL_h(\imat x),
\]
where $x\in\fg$ is the unique vector satisfying $d_eL_h(x)=X(h)$. On the other hand we have
\begin{align*}
 ((L_g)^*J)(X(h))&=d_{gh}L_{g^{-1}}\circ(J_{gh})\circ (d_hL_g(X(h)))\\
 &=d_{gh}L_{g^{-1}}\circ d_eL_{gh}(\imat x)=d_eL_h(\imat x)\\
 &=J_h(X(h)).
\end{align*}
\end{proof}


\begin{lemma}
The (almost) complex structure def{}ined above is integrable.
\end{lemma}

\begin{proof}
 This is trivial. Since $J$ is left invariant, it is enough to work in the Lie algebra $\fg$. But on $\fg$ the 
 complex structure is multiplication by $\imat$, hence the Nijenhuis tensor
\begin{align*}
 N_J(X,Y) &= [X,Y]+J[JX,Y]+J[X,JY]-[JX,JY] \\ 
&= [X,Y]+\imat [\imat X,Y]+\imat [X,\imat Y]-[\imat X,\imat Y]=0, 
\end{align*}
 for $X,Y\in\fg$.
\end{proof}

\begin{lemma}
 The two complex structures on $G$ coincide.
\end{lemma}
\begin{proof}
 It is enough to notice that the left translations are holomorphic maps, thus their dif{}ferential is complex linear. 
 Let 
 $g\in G$ be a point and $X$ a left invariant vector f{}ield on $G$, such that $X(g)=d_eL_g(x)$, $x\in\fg$. Then
 \[
 \imat X(g)=\imat d_eL_g(x)=d_eL_g(\imat x)=J_g(X(g)).
 \]
\end{proof}

So far we know that the natural complex structure $J$ on the Lie group $G=H_{\bC}\times\bC$ is 
left invariant and it is multiplication by $\imat$ on each tangent space. As above, let $G_\Gamma\subset G$ be
the subgroup of matrices whose elements belong to the lattice $\Gamma=\{a+b\zeta \ | \ \zeta=
e^{2\pi \imat /3}, a,b\in \bZ \}\subset\bC$. Since $J$ is left invariant, it def{}ines a complex structure on the quotient
$M=G_\Gamma\ba G$, which we denote again by $J$. Hence $(M,J)$ is a complex nilmanifold.

Next we show that $J$ is compatible with the $\bZ_3$-action def{}ined by (\ref{Z3-action}). The complex structure $J$ on 
$M$ is multiplication by $\imat$ at each tangent space $T_pM$, $p\in M$, since it comes from 
the complex structure on $G$. Let $\varphi\colon M\to M$ denote the $\bZ_3$-action, and consider the map
\[
d_p\varphi\colon T_pM\to T_{\varphi(p)}M.
\]
We claim that the map $\varphi$ can be lifted to a holomorphic action $\tilde{\varphi}$ of $\bZ_3$ on $G$. 
By taking global coordinates $(u_1,u_2,u_3,u_4)$ on $G$, $\tilde{\varphi}$ sends the generator
$\zeta\in\bZ_3$ to the diagonal matrix $\mathrm{diag}(\zeta,\zeta,\zeta^2,\zeta)$. Since $\tilde{\varphi}$ is linear, it coincides with its 
dif{}ferential $d_g\tilde{\varphi}\colon T_gG\to T_{\tilde{\varphi}(g)}G$. This is clearly a complex
linear map, i.e.
\begin{equation}\label{eq:456}
 d_g\tilde{\varphi}\circ J_g=J_{\tilde{\varphi}(g)}\circ d_g\tilde{\varphi}.
\end{equation}
This proves the claim. Since the complex structure $J$ on $M$ is multiplication 
by $\imat$ on each tangent space, (\ref{eq:456}) shows that we can write
\[
 d_p\varphi\circ J_p=J_{\varphi(p)}\circ d_p\varphi,
\]
showing that the complex structure commutes with the $\bZ_3$-action, hence descends to the quotient 
$\widehat{M}=M/\bZ_3$. We denote by $\widehat{J}$ the complex structure on $\widehat{M}$. 
Thus 
we have proved:

\begin{proposition}
 Let $M=G_\Gamma\ba G$ be as above and denote by $J$ the natural complex structure on $M$. 
 Then $(\widehat{M},\widehat{J})$ is a complex orbifold.
\end{proposition}

\begin{remark}
The complex nilmanifold $M$ is an example of an $8$-dimensional non-simply connected complex, symplectic and 
non-K\"ahler manifold, the symplectic form being given by \eqref{simplform}. Indeed, $M$ is
non-formal, hence it can not be K\"ahler. One can show that $(\widehat{M},\widehat{J},\widehat{\omega})$ is simply connected. Therefore we have an example of an 8-dimensional simply connected complex
and symplectic orbifold which is not K\"ahler. Indeed, one can see that $\widehat{M}$ is not formal \cite{FM}, while K\"ahler orbifolds are formal \cite{BFMT}.
\end{remark}

\begin{proposition} \label{prop:desingular}
 There exists a smooth complex manifold $(\overline{M},\overline{J})$ which is biholomorphic to $(\widehat{M},\widehat{J})$ outside a neighborhood of a singular point.
\end{proposition}

\begin{proof}
Let $p\in M$ be a f{}ixed point of the $\bZ_3$-action. Translating with an element $g\in G$, we 
can suppose that $p=(0,0,0,0)$ in our coordinates. Let $U\subset M$ be a neighborhood of $p$ and let
$\phi\colon U\to B$ be a holomorphic local 
chart, given by the exponential map (by holomorphic we mean with respect to the complex structure $J$). 
Here $B=B_{\bC^4}(0,\varepsilon)\subset\bC^4$. In these coordinates, the action of $\bZ_3$ can
be written as
$$
(u_1,u_2,u_3,u_4)\mapsto(\zeta u_1,\zeta u_2,\zeta^2 u_3,\zeta u_4).
$$
The local model for the singularity is thus $B/\bZ_3$. From now on, the desingularization process is 
formally analogous to that in \cite{FM}. We blow up $B$ at $p$ to obtain $\widetilde{B}$. The point
$p$ is replaced with a 
complex projective space $F=\bP^3=\bP(T_pB)$ on which $\bZ_3$ acts by
$$
[u_1:u_2:u_3:u_4]\mapsto[\zeta u_1:\zeta u_2:\zeta^2 u_3:\zeta u_4]=[u_1:u_2:\zeta u_3:u_4].
$$
Thus $\bZ_3$ acts on the exceptional divisor $F$ with f{}ixed locus $\{q\}\cup H$ where $q=[0:0:1:0]$ 
and $H=\{u_3=0\}$. Then one blows up $\widetilde{B}$ at $q$ and $H$ to obtain
$\widetilde{\widetilde{B}}$. The point $q$ is replaced by a 
projective space $H_1\cong\bP^3$. The normal bundle to $H\subset F\subset\widetilde{B}$ is the sum of the 
normal bundle of $H$ in $\bP^3$, which is $\mathcal{O}_{\bP^2}(1)$, and the restriction to $H$
of the normal bundle of $F$ in 
$\widetilde{B}$, which is $\mathcal{O}_{\bP^2}(-1)$. Hence the second blow up replaces the projective plane $H$ with a $\bP^1$-bundle over $\bP^2$ def{}ined as
$H_2=\bP(\mathcal{O}_{\bP^2}(1)\oplus\mathcal{O}_{\bP^2}(-1))$. 
The strict transform of $F\subset\widetilde{B}$ under the second blow up is the blow up $\widetilde{F}$ of $F$ at $q$, which is a $\bP^1$-bundle over $\bP^2=H$, actually 
$\widetilde{F}=\bP(\mathcal{O}_{\bP^2}\oplus\mathcal{O}_{\bP^2}(1))$. The resulting situation is depicted in Figure \ref{fig:2}. 

\begin{figure}[h!]\label{fig:2}
\begin{center}
    \includegraphics[width=9cm]{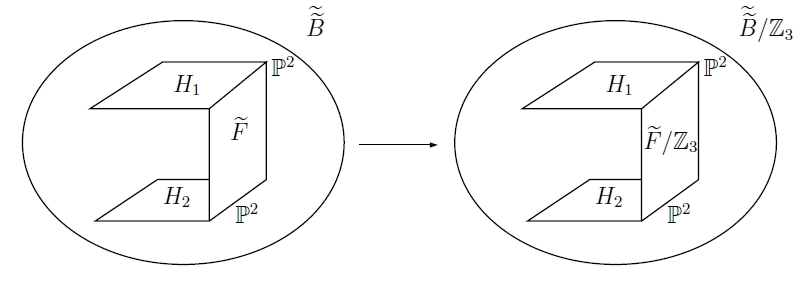}
\caption{The second blow-up and the $\bZ_3$-action}
\end{center}
\end{figure}

The f{}ixed point locus of the $\bZ_3$-action on $\widetilde{\widetilde{B}}$ consists of the 
two disjoint divisors $H_1$ and $H_2$. According to \cite{BPVdV}, page 82, the quotient
$\widetilde{\widetilde{B}}/\bZ_3$ is a smooth K\"ahler 
manifold. This provides a complex resolution of the singularity $B/\bZ_3$. Notice that the 
blowing up is performed with respect to the natural complex structure inherited from the ambient space. 
By resolving every singular point, we obtain a smooth complex manifold $(\overline{M},\overline{J})$.
\end{proof}

\begin{proposition}
 The complex manifold $(\overline{M},\overline{J})$ is simply connected.
\end{proposition}

\begin{proof}
 The proof is analogous to that of \cite[Proposition 2.3]{FM}.
\end{proof}

The desingularization process of \propref{prop:desingular} is completely similar to the symplectic resolution of \cite[Proposition 2.1]{FM}. 
However, the two blow ups are performed with respect to \textit{dif{}ferent} complex structures.
In the complex resolution, one uses the natural complex structure $\widehat{J}$ of $\widehat{M}$. In the 
symplectic resolution one uses a (local) complex structure obtained by using a K\"{a}hler model for
a neighborhood of a f{}ixed point which is not holomorphically equivalent to a local holomorphic chart 
for $\widehat{J}$. Indeed, this K\"ahler model is obtained by performing the following change of
coordinates in a small neighborhood of a f{}ixed point of the action:
\begin{equation}\label{cv1}
\left\{
\begin{array}{rcl}
w_1 & = & u_1\\
w_2 & = & \frac{1}{\sqrt{2}}(u_2+\imat\bar{u}_3)\\
w_3 & = & \frac{1}{\sqrt{2}}(\imat\bar{u}_2-u_3)\\
w_4 & = & u_4\\
\end{array}\right.
\end{equation}
Certainly, this is not holomorphic with respect to the natural complex structure $\widehat{J}$ on $\widehat{M}$.

Locally, we have the following situation: on a small neighborhood $\mathcal{U}$ of $0\in\bC^4$ (which is a f{}ixed point of the $\bZ_3$-action in suitable coordinates)
we have two complex structures, $J_1$ and $J_2$. The two complex structures are dif{}ferent, because the change of variables which brings one to the other is not holomorphic. As a consequence, the two
blow ups are dif{}ferent. In fact, the natural map that one would construct from one resolution to the other would not be even continuous. This becomes particularly clear when the blow up is
interpreted
as a symplectic cut, following Lerman and McDuf{}f (see for instance \cite{L}). The blow up of $\bC^n$ at $0$ can be thought of as removing a small ball of radius $\varepsilon$ centered at the origin
and then collapsing the f{}ibers of the Hopf f{}ibration in the boundary of the remaining set. But the f{}ibers of the Hopf f{}ibration (i.e.\ the intersections of the boundary of the ball, which is a
$S^{2n-1}$, with the ``complex'' lines through the origin) depend heavily on the complex structure of the ball.

\section{Proof of the main Theorem}\label{Section:Diffeomorphism}

In this section we prove that the smooth manifolds which underly the two resolutions $\overline{M}$ and $\widetilde{M}$ are dif{}feomorphic. This completes the proof of \thmref{main}.

\begin{proposition}
The symplectic and the complex resolution of the orbifold $(\widehat{M},\widehat{J},\widehat{\omega})$ are dif{}feomorphic.
\end{proposition}

\begin{proof} We work locally, in a small neighborhood of each f{}ixed point. We construct a smooth 
map which is the identity outside this small neighborhood and that does the right job inside the
neighborhood. The local model is 
thus a small ball $B_{\bC^4}(0,\delta)\subset\bC^4$ endowed with two dif{}ferent complex structures 
$J_1$ and $J_2$. There is a map $\Theta\colon B_{\bC^4}(0,\delta)\to B_{\bC^4}(0,\delta)$ which interchanges
the two complex structures, namely
$$
\Theta^*J_1=J_2.
$$
Notice that $\Theta$ can be composed with biholomorphisms on the right and on the left, thus is not unique. 
If we take $J_1$ as the complex structure on the ball induced by the natural complex
structure on $\widehat{M}$ and $J_2$ 
to be the complex structure associated to the local K\"ahler model used for the symplectic resolution, then $\Theta$ 
is given by (\ref{cv1}).
We introduce real coordinates $u_k=x_k+iy_k$ and $w_k=s_k+it_k$, $k=1,2,3,4$. 
In such coordinates, (\ref{cv1}) is an automorphism of $\bR^8$ written as
\[
\left\{
\begin{array}{rcl}
s_1 & = & x_1\\
t_1 & = & y_1\\
s_2 & = & \frac{1}{\sqrt{2}}(x_2+y_3)\\
t_2 & = & \frac{1}{\sqrt{2}}(y_2+x_3)\\
s_3 & = & \frac{1}{\sqrt{2}}(y_2-x_3)\\
t_3 & = & \frac{1}{\sqrt{2}}(x_2-y_3)\\
s_4 & = & x_4\\
t_4 & = & y_4
\end{array}\right.
\]

The associated matrix is
$$
\Theta=\begin{pmatrix}
 1 & 0 & 0 & 0 & 0 & 0 & 0 & 0\\
 0 & 1 & 0 & 0 & 0 & 0 & 0 & 0\\
 0 & 0 & \frac{1}{\sqrt{2}} & 0 & 0 & \frac{1}{\sqrt{2}} & 0 & 0\\
 0 & 0 & 0 & \frac{1}{\sqrt{2}} & \frac{1}{\sqrt{2}} & 0 & 0 & 0\\
 0 & 0 & 0 & \frac{1}{\sqrt{2}} & -\frac{1}{\sqrt{2}} & 0 & 0 & 0\\
 0 & 0 & \frac{1}{\sqrt{2}} & 0 & 0 & -\frac{1}{\sqrt{2}} & 0 & 0\\
 0 & 0 & 0 & 0 & 0 & 0 & 1 & 0\\
 0 & 0 & 0 & 0 & 0 & 0 & 0 & 1
\end{pmatrix}.
$$
The matrix $\Theta$ belongs to $\mathrm{SO}(8)$. To construct the dif{}feomorphism we 
will f{}ind an isotopy $\{\Theta_t\}_{t\in[0,1]}$, such that $\Theta_0$ is the identity
$\mathrm{Id}\in\mathrm{SO}(8)$ and $\Theta_1=\Theta$. 
In this way we get a path of complex structures $J_{t+1}=\Theta_t^*J_1$ connecting $J_1$ and $J_2$.
To do this we must produce a smooth path in $\mathrm{SO}(8)$ between the identity matrix and $\Theta$, 
which is furthermore equivariant with respect to the $\bZ_3$-action. In fact it is enough to f{}ind a smooth
$\bZ_3$-equivariant path in $\mathrm{SO}(4)$ connecting the identity to the matrix
$$
\theta=\begin{pmatrix}
 \frac{1}{\sqrt{2}} & 0 & 0 & \frac{1}{\sqrt{2}}\\
 0 & \frac{1}{\sqrt{2}} & \frac{1}{\sqrt{2}} & 0\\
 0 & \frac{1}{\sqrt{2}} & -\frac{1}{\sqrt{2}} & 0\\
 \frac{1}{\sqrt{2}} & 0 & 0 & -\frac{1}{\sqrt{2}}
\end{pmatrix}
$$
In the coordinates $(s_2,t_2,s_3,t_3)$ spanning the $\bR^4$ of interest, the $\bZ_3$-action can be written as
$$
\Upsilon=\begin{pmatrix}
 -\frac{1}{2} & -\frac{\sqrt{3}}{2} & 0 & 0\\
 \frac{\sqrt{3}}{2} & -\frac{1}{2} & 0 & 0\\
 0 & 0 & -\frac{1}{2} & \frac{\sqrt{3}}{2}\\
 0 & 0 & -\frac{\sqrt{3}}{2} & -\frac{1}{2}
\end{pmatrix}
$$
under the natural inclusion $\mathrm{U}(2)\hookrightarrow\mathrm{SO}(4)$. We must 
ensure that the path $\{\Theta_s\}\subset\mathrm{SO}(4)$ satisf{}ies $\Theta_s\circ\Upsilon=\Upsilon\circ\Theta_s$, 
for every $s\in[0,1]$. We do this explicitly. First notice that $\theta=P\theta'$, where
$$
P=\begin{pmatrix}
 0 & 0 & 0 & 1\\
 0 & 0 & 1 & 0\\
 0 & 1 & 0 & 0\\
 1 & 0 & 0 & 0
\end{pmatrix}
\quad \textrm{and}\quad
\theta'=\begin{pmatrix}
 \frac{1}{\sqrt{2}} & 0 & 0 & -\frac{1}{\sqrt{2}}\\
 0 & \frac{1}{\sqrt{2}} & -\frac{1}{\sqrt{2}} & 0\\
 0 & \frac{1}{\sqrt{2}} & \frac{1}{\sqrt{2}} & 0\\
 \frac{1}{\sqrt{2}} & 0 & 0 & \frac{1}{\sqrt{2}}
\end{pmatrix}.
$$
The matrix $\theta'$ is the image, under the exponential map $\exp\colon\mathfrak{so}(4)\to \mathrm{SO}(4)$, of the matrix $\frac{\pi}{4}Q$, where
$$
Q=\begin{pmatrix}
 0 & 0 & 0 & -1\\
 0 & 0 & -1 & 0\\
 0 & 1 & 0 & 0\\
 1 & 0 & 0 & 0
\end{pmatrix}.
$$
Thus a smooth path in $\mathrm{SO}(4)$ between the identity and $\theta'$ is given by the image of the straight line in $\mathfrak{so}(4)$ joining the zero matrix with $Q$,
$$
\begin{array}{rclcc}
\gamma\colon & [0,\pi/4] & \rightarrow & \mathrm{SO}(4)\\
& s & \mapsto & \textrm{exp}(sQ) 
\end{array}
$$
One sees that, for every $s\in[0,\pi/4]$, $\gamma(s)\circ\Upsilon=\Upsilon\circ\gamma(s)$, hence $\gamma(s)$ 
is $\bZ_3$-equivariant. Now consider the matrix $P$. We juxtapose the following three paths
in order to join $P$ with the 
identity matrix:
$$
P_1(s)=\begin{pmatrix}
 0 & 0 & \sin(\pi s/2) & \cos(\pi s/2)\\
 0 & 0 & \cos(\pi s/2) & -\sin(\pi s/2)\\
 \sin(\pi s/2) & \cos(\pi s/2) & 0 & 0\\
 \cos(\pi s/2) & -\sin(\pi s/2) & 0 & 0
\end{pmatrix},
$$
$$
P_2(s)=\begin{pmatrix}
 \sin(\pi s/2) & 0 & \cos(\pi s/2) & 0\\
 0 & \sin(\pi s/2) & 0 & -\cos(\pi s/2)\\
 \cos(\pi s/2) & 0 & -\sin(\pi s/2) & 0\\
 0 & -\cos(\pi s/2) & 0 & -\sin(\pi s/2)
\end{pmatrix},
$$
$$
P_3(t)=\begin{pmatrix}
 1 & 0 & 0 & 0\\
 0 & 1 & 0 & 0\\
 0 & 0 & -\cos(\pi s) & \sin(\pi s)\\
 0 & 0 & -\sin(\pi s) & -\cos(\pi s)
\end{pmatrix}.
$$
Again, a computation shows that $P_i(s)\circ\Upsilon=\Upsilon\circ P_i(s), \ \forall s\in[0,1]$, $i=1,2,3$. 
Hence the path $P(s)=P_1\ast P_2\ast P_3(s)$ satisf{}ies $P(0)=P$, $P(1)=\mathrm{Id}$ and is
$\bZ_3$-equivariant. The path 
$\theta(s)=P(1-s)\theta'$ satisf{}ies $\theta(0)=\theta'$ and $\theta(1)=\theta$. Finally the path $\Psi=\gamma \ast \theta$ 
satisf{}ies $\Psi(0)=\mathrm{Id}$ and $\Psi(1)=\theta$. However $\Psi$ is not globally smooth, because at 
the concatenation points it is only continuous. To smooth it, we proceed as follows. Let $0<s_1<\ldots<s_{n-1}<s_n<1$ denote the points 
in which the resulting path has a cusp. Consider a smooth, increasing function $h\colon [0,1]\to[0,1]$
such that there exist intervals 
$\mathcal{J}_i=(t_i-\varepsilon,t_i+\varepsilon)$, $0<t_1<\ldots<t_{n-1}<t_n<1$ with $h(t)=s_i$ for 
$t\in\mathcal{J}_i$. Def{}ine a new path $\Theta_t=\Psi_{h(t)}$. Clearly $\Psi$ and $\Theta$ have the
same image. Then $\Theta_t$ is
a smooth, $\bZ_3$-equivariant path in $\mathrm{SO}(4)$ connecting $\theta$ with the identity matrix. Viewing 
it as a path in $\mathrm{SO}(8)$ we obtain the isotopy $\Theta_t$ such that $\Theta_0=\mathrm{Id}$ 
and $\Theta_1=\Theta$. Thus $\Theta_0^*J_1=J_1$ and $\Theta_1^*J_1=J_2$. 
We also endow the ball with the standard metric. Since $\bZ_3\subset \mathrm{SO}(8)$, $\bZ_3$ acts by isometries.

We are ready to def{}ine the dif{}feomorphism between the two resolutions. Notice that the 
expression of the $\bZ_3$-action is the same in the two sets of coordinates $(u_1,\ldots,u_4)$ and
$(w_1,\ldots,w_4)$. Thus when we blow up we get, 
in both cases, an exceptional divisor $\bP^3$ with one f{}ixed point $q=[0:0:1:0]$ and one f{}ixed hyperplane 
$H=\{u_3=0\}=\{w_3=0\}$. The dif{}ferential of $\Theta$ at $0\in B_{\bC^4}(0,\delta)$, which we
denote $d_0\Theta$, def{}ines an 
automorphism of the exceptional divisor (when we projectivize the action), which f{}ixes $q$ and maps $H$ to itself 
($d_0\Theta$ is $(J_1,J_2)$-holomorphic, meaning that $d_0\Theta\circ J_1=J_2\circ
d_0\Theta$). Thus $d_0\Theta$ also 
lifts to the second blow-up, hence to a map between the two exceptional divisors. Let $\rho\colon\bR\to[0,1]$ be a 
cut-of{}f function, i.e.\ a $C^{\infty}$ function, which is identically $0$ on
$(-\infty,0]$ and identically $1$ on 
$[1,\infty)$. Using the metric on the ball, the dif{}feomorphism $f$ can then be def{}ined as 
follows:
$$
f(x)=\left\{
\begin{array}{ll}
x & \textrm{if} \ |x|>\frac{2\delta}{3}\\
& \\
\Theta_t(x) & \textrm{if} \ \frac{\delta}{3}<|x|<\frac{2\delta}{3}\\
& \\
\Theta(x) & \textrm{if} \ |x|<\frac{\delta}{3}\\
\end{array}
\right.
$$
where $t=\rho\big(\big(\frac{2\delta}{3}-|x|\big)\frac{3}{\delta}\big)$.
By what we have said, $f\colon\widehat{M}\to \widehat{M}$ lifts to a dif{}feomorphism $\tilde f\colon \overline{M}\to \widetilde{M}$.
\end{proof}

\begin{corollary}
The manifold $\widetilde{M}$ is a simply connected, $8$-dimensional,  non-formal manifold that admits both
complex and symplectic structures, but which carries no K\"ahler metric.
\end{corollary}

\end{document}